\documentclass[useAMS,usenatbib,12pt,epsfig]{article}
\usepackage[X2,T2A]{fontenc}
\usepackage[english]{babel}
\usepackage{amssymb}
\usepackage{graphicx}
\usepackage[dvips]{color}

\begin{document}
\selectlanguage{english}

\thispagestyle{empty}

\bigskip
\bigskip
\bigskip
\begin{center}
\begin{Large}
{An integral representation of divisor function. 
\\An equation for prime numbers.}
\end{Large}
\end{center}
\bigskip
\bigskip
\bigskip
\begin{center}
{Kholupenko E.E.}
\end{center}
Ioffe Physical-Technical Institute, St.-Petersburg 194021, Russia
\bigskip
\bigskip
\bigskip
\bigskip
\begin{center}
\bf{Abstract}
\end{center}
{
A representation of divisor function $\tau(n)\equiv \sigma_{0}(n)$ by means of  
logarithmic residue of a function of complex variable is suggested. 
This representation may be useful theoretical instrument for further investigations 
of properties of natural numbers.
}
\\ 
\bigskip
\\
{Keywords: prime numbers, divisor function, trigonometric sums, analytic number theory}

\newpage
\section{Introduction}
Properties of natural numbers are the fundamental question of mathematics. 
Investigations of these properties in the Ancient Times (e.g. the infinitude 
of prime numbers which was shown by Euclid, Sieve of Eratosthenes for prime 
numbers, consideration of Diophantine equations) 
led to the birth of elementary theory of numbers. 
Later investigations (e.g. Fermat's Last theorem, investigations of 
zeta-function by Euler and Riemann) have led to the development 
the great branch of mathematics - the analytical theory of numbers. 
The important problems of modern number theory are the problems 
of effective search of prime numbers and complexity (and factorization) or 
primality of arbitrary natural number. 
Today essential efforts are applied for consideration of these questions: 
many fast methods for factorization (e.g. Pollard $\rho$-algorithm \cite{b1}, 
elliptic curve method \cite{b2}, quadratic sieve algorithm \cite{b3}, general number field sieve \cite{b4}) 
and tests of primality (e.g. Miller-Rabin test \cite{b5}, Lucas-Lehmer test \cite{b6}, AKS test \cite{b7}) 
are developed, projects on search of large prime numbers are running 
(GIMP \cite{b8}, PrimeGride \cite{b9}, Seventeen or Bust \cite{b10} and others), 
databases of prime numbers are updated (e.g. The List of Largest Known Primes \cite{b11}, OEIS \cite{b12}). 
Such interest is due to not only desire to understand the nature of 
numbers and to resolve some fundamental problems of mathematics but also 
due to opportunity of practical application of properties of natural numbers, 
e.g. for the cryptographic systems with public key (for example RSA \cite{b13}), 
hash-tables, and generation of pseudorandom numbers (e.g. Mersenne twister \cite{b14}).

The aim of this paper is to suggest a theoretical instrument for 
investigation of quantity of divisors of natural number, which may be 
useful for further investigations of properties of natural numbers.

\section{Integral representation of divisor function}
\begin{figure}
\begin{center}
\includegraphics[width=16cm]{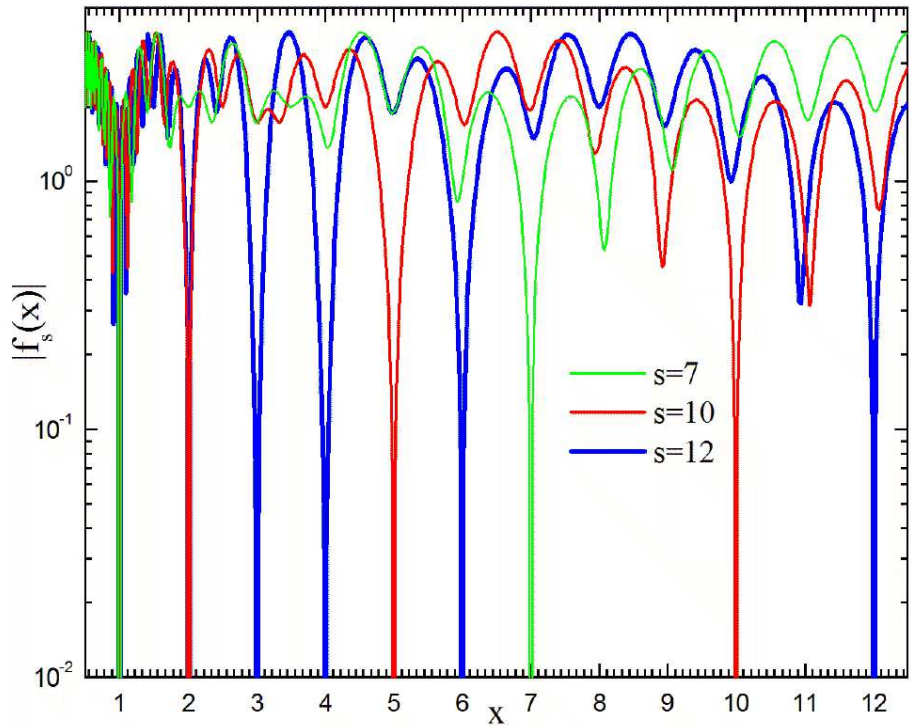}
\end{center}
\caption{Dependence of module of function $f_{s}(x)$ on real variable $x$.}
\label{fig1}
\end{figure}

Main idea is similar to the one of method of trigonometric sums \cite{b15}. 
Using the fact that $\exp\left(2\pi i x\right)$ equals unity only for 
integer value of real number $x$ one may say that function: 
\begin{equation}
f_{s}(x)=\exp\left(2\pi i x\right)+\exp\left(2\pi i s/x\right) - 2
\end{equation}
has zeros only for values of $x$ which are divisors of natural number $s$ (we will 
denote these numbers as $x^{s}_{j}$, where superscript $s$ highlights dependence 
on natural number $s$,
subscript $j$ denotes integer sequence number of corresponding zero). 
Note that 
\begin{equation}
\exp\left(2\pi i x^{s}_{j}\right)=1,~~~\exp\left(2\pi i s/x^{s}_{j}\right)=1
\label{x_s_condition}
\end{equation}
Unfortunately, direct use of $f_{s}(x)$ (e.g. search of global minima of $f_{s}(x)$ by using numerical iteration methods) 
for factorization of natural number $s$ or for counting number of divisors of $s$ is almost impossible because of frequent 
oscillations of function $f_{s}(x)$ (see fig. \ref{fig1}). 

Number of zeros of this function in the interval $[1;s]$ gives us the 
number of divisors of $s$. To count the number of zeros of $f_{s}$ one 
can use logarithmic residue of $f_{s}$. In general case the logarithmic residue is the following (e.g. \cite{b16}):
\begin{equation}
\gamma \left[f, L\right] = {1\over 2\pi i}\int_{L} 
{f'\left(z\right) dz \over f\left(z\right)}=N_{f}-P_{f}
\label{gamma_f_L_def}
\end{equation}
where $f$ is a function of complex number $z$, $L$ is a closed 
contour on the complex plane $\mathbb{C}$, $N_{f}$ is the number of 
zeros of function $f$ within the contour $L$ (taking into account the orders of zeros), 
$P_{f}$ is the number of poles of function $f$ within the contour $L$ (taking into account the orders of poles). 
If function $f$ has no poles within contour $L$ (e.g. $f$ is the smooth functionof $z$) then one may write:
\begin{equation}
\gamma \left[f, L\right] = \sum_{z_{j} \in \mathbb{L}: f\left(z_{j}\right)=0}k_{i}
\label{gamma_f_L}
\end{equation}
where $\mathbb{L}$ is the field within the contour $L$ (i.e. $\mathbb{L}$ is the field limited by the contour $L$ 
on the complex plane $\mathbb{C}$), $z_{j}$ are zeros of function $f$ in the field $\mathbb{L}$, 
$k_{j}$ is the order of zero $z_{j}$, i.e. in the vicinity of $z_{j}$ the function $f\left(z\right)$ can be presented as 
\begin{equation}
f\left(z\right)={f^{\left(k_{j}\right)}\left(z_{j}\right)\over k_{j}!}\left(z-z_{j}\right)^{k_{j}}+O\left(\left(z-z_{j}\right)^{k_{j}+1}\right)
\end{equation}

\begin{figure}
\begin{center}
\includegraphics[width=12cm]{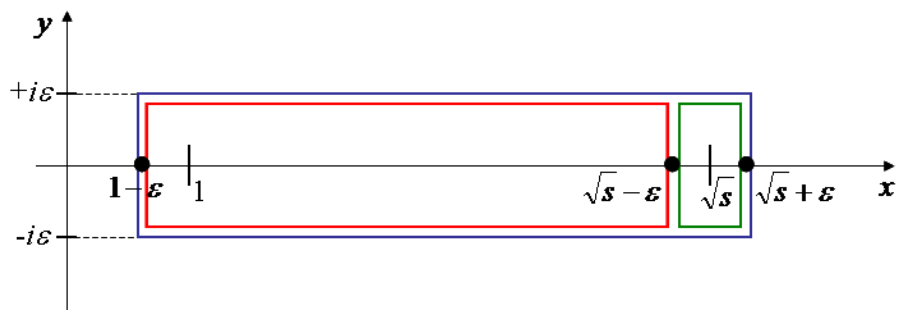}
\end{center}
\caption{Examples of contours $L$ in complex plane for calculation of logarithmic residue of function $f_{s}(x)$.}
\label{fig2}
\end{figure}

Using (\ref{gamma_f_L}) one can construct the expression for counting number 
of zeros of function $f_{s}\left(z\right)$ in the interval $[1;\sqrt{s}]$. 
Let us define the following functions:
\begin{equation}
\tau_{-} (s)=\gamma \left[f_{s}, L^{-}_{s}\right]\;,~~~~
\tau_{c} (s)=\gamma \left[f_{s}, L^{c}_{s}\right]\;,~~~~
\tau_{+} (s)=\gamma \left[f_{s}, L^{+}_{s}\right]
\label{taus}
\end{equation}
where $L^{-}_{s}$ is the contour given (e.g.) by the following points $z$ (see fig. \ref{fig2}):
\begin{eqnarray}
z=1-\epsilon+iy,~~y\in \left[-\epsilon;+\epsilon\right]\nonumber\\
x+i\epsilon,~~x\in \left[\left(1-\epsilon\right);\left(\sqrt{s}-\epsilon\right)\right]\nonumber\\
\left(\sqrt{s}-\epsilon\right)+iy,~~y\in \left[+\epsilon;-\epsilon\right]\nonumber\\
x-i\epsilon,~~x\in \left[\left(\sqrt{s}-\epsilon\right);\left(1-\epsilon\right)\right]
\end{eqnarray}
$L^{c}_{s}$ is the contour given (e.g.) by the following points $z$ (see fig. \ref{fig2}):
\begin{eqnarray}
z=\sqrt{s}-\epsilon+iy,~~y\in \left[-\epsilon;+\epsilon\right]\nonumber\\
x+i\epsilon,~~x\in \left[\left(\sqrt{s}-\epsilon\right);\left(\sqrt{s}+\epsilon\right)\right]\nonumber\\
\left(\sqrt{s}+\epsilon\right)+iy,~~y\in \left[+\epsilon;-\epsilon\right]\nonumber\\
x-i\epsilon,~~x\in \left[\left(\sqrt{s}+\epsilon\right);\left(\sqrt{s}-\epsilon\right)\right]
\end{eqnarray}
$L^{+}_{s}$ is the contour given by the following points $z$ (see fig. \ref{fig2}):
\begin{eqnarray}
z=1-\epsilon+iy,~~y\in \left[-\epsilon;+\epsilon\right]\nonumber\\
x+i\epsilon,~~x\in \left[\left(1-\epsilon\right);\left(\sqrt{s}+\epsilon\right)\right]\nonumber\\
\left(\sqrt{s}+\epsilon\right)+iy,~~y\in \left[+\epsilon;-\epsilon\right]\nonumber\\
x-i\epsilon,~~x\in \left[\left(\sqrt{s}+\epsilon\right);\left(1-\epsilon\right)\right]
\end{eqnarray}
where $\epsilon$ is the positive real number. Parameter $\epsilon$ should be small enough to exclude 
zeros $z^{s}_{j}$ ($z^{s}_{j}$ denote zeros of $f_{s}(z)$ in complex plane)  which do not lie on the real axis. 
Note that $\tau_{+}(s)=\tau_{-}(s)+\tau_{c}(s)$.

Studying values of $\tau_{-}(s)$, $\tau_{c}(s)$ and $\tau_{+}(s)$ one should calculate the orders 
of zeros $z^{s}_{j}$ which lie on the real axis (i.e. $z^{s}_{j}=x^{s}_{j}$). 
This demands calculation of derivative of $f_{s}\left(z\right)$:
\begin{equation}
f'_{s}\left(z\right)=2\pi i\left[\exp\left({2\pi i z}\right)-{s \over z^2}\exp\left({2\pi i s\over z}\right)\right]
\end{equation}
Taking into account (\ref{x_s_condition}) one may write
\begin{equation}
f'_{s}\left(x^{s}_{j}\right)=2\pi i\left[1-{s \over \left(x^{s}_{j}\right)^2}\right]
\label{f1_simplified}
\end{equation}
Thus $f'_{s}\left(x^{s}_{j}\right)\neq 0$ for $x^{s}_{j}$ $\in$ $\left[1;\sqrt{s}\right)$, and therefore 
all zeros $x^{s}_{j}$ in interval $\left[1;\sqrt{s}\right)$ have order which equals unity. If $\sqrt{s}$ is 
also divisor of $s$ then one should investigate the order of this zero. From (\ref{f1_simplified}) it is 
evident that this order is not equal unity. Thus one should find the second order of $f_{s}(z)$:
\begin{equation}
f''_{s}\left(z\right)=2\pi i \left[2\pi i\exp\left({2\pi i z}\right)+{2s \over z^3}\exp\left({2\pi i s\over z}\right)+2\pi i{s^2 \over z^4}\exp\left({2\pi i s\over z}\right)\right]
\label{f2}
\end{equation}
If the value $\sqrt{s}$ is an integer number (this case is interesting for us) then one may write:
\begin{equation}
f''_{s}\left(\sqrt{s}\right)=-8\pi ^2 + 4\pi i s^{-1/2} \neq 0
\label{f2_calculated}
\end{equation}
Expression (\ref{f2_calculated}) means that the order of zero in point $\sqrt{s}$ equals 2 (if $\sqrt{s}$ is integer and divisor of number $s$) 
or 0 (if $\sqrt{s}$ is not integer and not a zero of function $f_{s}\left(z\right)$), i.e. $\tau_{c}(s)=0$ or 2.

Knowledge of behaviour of functions $\tau_{-}(s)$, $\tau_{c}(s)$ and $\tau_{+}(s)$ allows us to calculate the number of 
divisors $\tau\left(s\right)$ of integer number $s$. Taking into account symmetry of quantity of divisors of $s$ relative to the number $\sqrt{s}$ 
one may write: 
\begin{equation}
\tau\left(s\right) = 2\tau_{-}(s)+{1\over 2}\tau_{c}(s) = {3 \over 2}\tau_{-}(s)+{1\over 2}\tau_{+}(s) = 2\tau_{+}(s)-{3\over 2}\tau_{c}(s)
\label{tau_representation}
\end{equation}

Explicit expressions for $\tau_{-}(s)$, $\tau_{c}(s)$ and $\tau_{+}(s)$ are the following:
\begin{equation}
\tau_{-} (s)=\int_{L^{-}_{s}}
{\exp\left({2\pi i z}\right)-\left({s / z^2}\right)\exp\left({2\pi i s/z}\right) \over 
\exp\left({2\pi i z}\right) + \exp\left({2\pi i s/z}\right)-2}dz
\label{tau_minus}
\end{equation}
\begin{equation}
\tau_{c} (s)=\int_{L^{c}_{s}}
{\exp\left({2\pi i z}\right)-\left({s / z^2}\right)\exp\left({2\pi i s/z}\right) \over 
\exp\left({2\pi i z}\right) + \exp\left({2\pi i s/z}\right)-2}dz
\label{tau_c}
\end{equation}
\begin{equation}
\tau_{+} (s)=\int_{L^{+}_{s}}
{\exp\left({2\pi i z}\right)-\left({s / z^2}\right)\exp\left({2\pi i s/z}\right) \over 
\exp\left({2\pi i z}\right) + \exp\left({2\pi i s/z}\right)-2}dz
\label{tau_plus}
\end{equation}

\section{Prime number equation}
One of the interesting cases of using representation (\ref{tau_representation}-\ref{tau_plus}) appears 
when $s$ is a prime number. In this case $\tau (s)=2$ and one can write the equation for prime numbers 
using equivalent condition $\tau_{+} (s)=1$. To prove the equivalence of $\tau (s)=2$ and $\tau_{+} (s)=1$ let us consider the following: 
\\1) From $\tau (s)\ge 2$ and $\tau (s)\le 2\tau_{+} (s)$ (from (\ref{tau_representation})) 
one can show that $\tau_{+} (s)=1 \Rightarrow \tau (s)=2$,
\\2) From $\tau_{+} (s) \ge 1$, $\tau_{+} (s)\ne 2$, and $\tau_{+} (s)\le \tau (s)$
one can show that $\tau (s)=2 \Rightarrow \tau_{+} (s)=1$.

Thus one can write the following equation for prime numbers $p$:
\begin{equation}
\int_{L^{+}_{p}}
{\exp\left({2\pi i z}\right)-\left({p / z^2}\right)\exp\left({2\pi i p/z}\right) \over 
\exp\left({2\pi i z}\right) + \exp\left({2\pi i p/z}\right)-2}dz=1
\label{prime_number_equation}
\end{equation}

\section{Conclusion}
Formulae (\ref{tau_representation}-\ref{tau_plus}) for an integral representation of divisor function 
and equation (\ref{prime_number_equation}) for prime numbers have been shown in the present paper. 
Unfortunately the calculation of expressions (\ref{tau_representation}-\ref{tau_plus}) and 
(\ref{prime_number_equation}) for the typical 
values of $s$ and $p$ (which are needed in practical problems) demands huge computational time 
and resources. Thus direct application of formulae 
(\ref{tau_representation}-\ref{tau_plus}) and (\ref{prime_number_equation}) is impossible. 
Nevertheless, these formulae may be useful for further investigations of properties of 
natural numbers and solving of practical and fundamental problems of number theory 
such as construction of new tests of primality and proof of the infinitude of twin primes.

\newpage

\newpage


\newpage
\begin{thebibliography}{}
\bibitem{b1}
Pollard J. M., BIT Numerical Mathematics 15 (3): 331 - 334 (1975)
\bibitem{b2}
Lenstra Jr. H. W., Annals of Mathematics 126 (2): 649 - 673 (1987)
\bibitem{b3}
Pomerance C., Computational Methods in Number Theory, Part I, Lenstra Jr. H.W. and Tijdeman R., eds., Math. Centre Tract 154, Amsterdam, pp 89-139 (1982)
\bibitem{b4}
Kleinjung T.,  Mathematics of Computation 75 (256): 2037 - 2047 (2006)
\bibitem{b5}
Rabin M. O., Journal of Number Theory 12 (1): 128 - 138 (1980)
\bibitem{b6}
Lehmer D. H., Annals of Mathematics, 31, 419 - 448 (1930) 
\bibitem{b7}
Agrawal M., Kayal N., Saxena N., Annals of Mathematics 160 (2): 781 - 793 (2004)
\bibitem{b8} 
Great Internet Mersenne Prime Search, http://www.mersenne.org/
\bibitem{b9} 
http://www.primegrid.com/
\bibitem{b10} 
http://www.seventeenorbust.com/
\bibitem{b11} 
Caldwell, C. K., The List of Largest Known Primes, http://primes.utm.edu/primes/
\bibitem{b12} 
The On-Line Encyclopedia of Integer Sequences, http://oeis.org/
\bibitem{b13}
Rivest R., Shamir A., Adleman L., Communications of the ACM 21 (2): 120 - 126 (1978)
\bibitem{b14}
Matsumoto M., Nishimura T., ACM Trans. on Modeling and Computer Simulations 8 (1): 3-30 (1998)
\bibitem{b15}
Vinogradov I.M.,   "The method of trigonometric sums in the theory of numbers" , Interscience  (1954)  (Translated from Russian)
\bibitem{b16}
Goncharov V. L., Teoriya funkcij kompleksnogo peremennogo ("The theory of functions of complex variable"), in russian, 
publishing company ''Prosveschenie'', Moscow (1955)
\end{thebibliography}
\end{document}